\documentclass[12pt]{amsart}
\input {epsf}
\usepackage{amssymb}
\usepackage{amsmath}
\usepackage{graphicx}         %%%%
\usepackage[all]{xy}
\newcommand\Z{{\mathbb{Z}}}

\newtheorem{theorem}{Theorem}[section]

\newtheorem{lemma}[theorem]{Lemma}

\begin{document}
%\date{\today}

\title{Yamada Polynomial and Khovanov Cohomology}

\author[Vershinin]{V.~V.~Vershinin}
\address{D\'epartement des Sciences Math\'ematiques,
Universit\'e Montpellier II, Place Eug\'ene Bataillon, 34095
Montpellier cedex 5, France} \email{vershini@math.univ-montp2.fr}
\address{ Sobolev Institute of Mathematics, Novosibirsk, 630090,
Russia} \email{versh@math.nsc.ru}
\author[Vesnin]{A.~Yu.~Vesnin}
\address{Sobolev Institute of Mathematics, Novosibirsk, 630090,
Russia} \email{vesnin@math.nsc.ru}

\begin{abstract}
For any graph $G$ we define bigraded cohomology groups whose graded Euler
characteristic is a multiple of the Yamada polynomial of $G$.
\end{abstract}



\thanks{%$^{*}$
The first named author was supported
 in part by the
 by CNRS-NSF grant No~17149, INTAS grant No~03-5-3251
and the ACI project ACI-NIM-2004-243 "Braids and Knots".}

\thanks{The second named author is supported by
INTAS grant ``CalcoMet-GT'' 03-51-3663 and by the grant of SB RAN}

\maketitle

\section{Introduction}

Mikhail Khovanov \cite{K00} constructed a bigraded homology group for
links such that its graded Euler characteristic is equal to the Jones polynomial.
The essential point of the construction
is the state sum formula for the Jones polynomial suggested in
\cite{Ka}. Since then many aspects of Khovanov's construction were
studied and generalized in various ways (see \cite{BN05, K05, V04}
and references therein). The existence of state sum descriptions for
diverse polynomial invariants gives the possibility to make
analogues of Khovanov's construction in other situations. In
particular, the similar constructions can be done in the cases of
some polynomial invariants of graphs.

In \cite{HR05} L.~Helme-Guizon and Y.~Rong constructed a cohomology
theory that categorifies the chromatic polynomial for graphs, i.~e.
the graded Euler characteristic of the constructed chain complex and
the corresponding homology groups is the chromatic polynomial.
E.F.~Jas\-so-Hernandez and Y.~Rong \cite{JR05} did the same for the
Tutte polynomial of graphs. It is natural to ask if similar
constructions can be made for other graph polynomials.

In the present paper we suggst a categorification for the two variables
Yamada
polynomial of graphs, which is universal among graph invariants
satisfying the deletion-contraction relation.  More precisely, for
each graph $G$ we define bigraded cohomology groups whose Euler
characteristic is a multiple  of the Yamada polynomial of $G$.

Our construction starts in Section~3 with rewriting the Yamada
polynomial in a state sum which is more friendly for a chain
complex set up. Section~4 is devoted to the construction of our
cohomology theory and proving main properties. In many aspects our
construction follows the ideas of the works \cite{HR05} and
\cite{JR05}. In Section~5 we give a simple example that illustrates
the construction.

\section{Polynomials of graphs \label{sec:pol_gr}}

Let $G$ be a finite graph with the vertex set $V(G)$ and the edge
set $E(G)$. For a given edge $e\in E(G)$ let $G-e$ be the graph
obtained from $G$ by {\em deleting} the edge $e$, and $G/e$ be the
graph obtained by {\em contracting} $e$ to a vertex (i.~e. by
deleting $e$ and identifying its ends to a single vertex). Recall
that $e$ is called a {\em loop} if $e$ joins a vertex to itself, and
is called an {\em isthmus} if its deleting from $G$ increases the
number of connected components of the graph. Two graphs are said to
be {\em 2-isomorphic} if there is a bijection between their edges
which induces a one-to-one correspondence between their cycles
\cite{Wh}. A graph function is said to be {\em $2$-invariant} if it
assigns to $2$-isomorphic graphs the same value.

Let $f$ be a $2$-invariant graph function with values in some ring
$R$. We will assume that a function $f$ satisfies the following
conditions:

$1^0$. {\em ``Deletion-contraction relation''.} If an edge $e$ is
not a loop or an isthmus then $f(G) = A  f(G / e) + B f(G - e)$,
where the coefficients $A\in R$ and $B\in R$ do not depend on
the choice of $e$.

$2^0$. If $H \cdot K$ is a union of two subgraphs $H$ and $K$ which
have only a common vertex then $f(H \cdot K) = C f(H) f(K)$, where
the coefficient $C\in R$ does not depend on the subgraphs $H$ and $K$.

$3^0$. If $T_1$ is a tree with a single edge on two vertices then
$f(T_1) = D$, for some $D \in R$.

$4^0$. If $L_1$ is a single-vertex graph with only loop then
$f(L_1)=E$, for some $E \in R$.

\smallskip

Thus, $f$ is determined by five coefficients $A$, $B$, $C$, $D$, and
$E$. Applying the above properties of the function $f$, we may
immediately calculate $f$ for
the simplest classes of graphs:\\
-- if $T_n$ is a connected tree with $n$ edges, then
$f(T_n) = C^{n-1} D^n$;\\
-- if $L_n$ is a single vertex with $n$ edges, then $f(L_n)
 = C^{n-1} E^n$; \\
-- if $D_n$ consists of two vertices joined by multiple $n$ edges,
then
$$
f(D_n) = B^{n-1} D + A E \frac{B^{n-1} - (C E)^{n-1} }{B - C E} ;
$$
-- if a graph $P_n$ is a simple cycle with $n$ edges, then
$$
f(P_n) = A^{n-1} E + B D \frac{A^{n-1} - (C D)^{n-1}}{A - C D} .
$$

\smallskip

For some particular values of coefficients (see the table below) the
graph function $f$ coincides with classical graph invariants: the
Tutte polynomial, the chromatic polynomial, and the flow polynomial
(see \cite{Tu} for definitions); the Negami polynomial \cite{Ne};
and the Yamada polynomial \cite{Ya}.

\begin{center}
\begin{tabular}[c]{|c|ccccc|}
\hline The polynomial & $A$ & $B$ & $C$ & $D$ & $E$
\\ \hline \hline Tutte polynomial $T(G;x,y)$ & $1$ & $1$ & $1$ & $x$ & $y$
\\ \hline
chromatic polynomial $P(G;\lambda)$ & $-1$ & $1$ & $1 \over \lambda$
& $\lambda(\lambda-1)$ & $0$
\\ \hline
flow polynomial $F(G;\lambda)$ & $1$ & $-1$ & $1$ & $0$ &
$\lambda-1$ \\ \hline Negami polynomial $N(G;t,x,y)$ & $x$ & $y$ &
$1 \over t$ & $t(x+ty)$ & $t(x+y)$ \\ \hline Yamada polynomial
$h(G;x,y)$ & $1$ & $- {1 \over x}$ & $1 \over x$ & $0$ & $xy-1$ \\
\hline
\end{tabular}
\end{center}

\smallskip

One can try to use well-known state sum formulae for these
polynomials to  categorify them. For each $S \subseteq E(G)$ let
$[G:S]$ be the graph whose vertex set is $V(G)$ and whose edge set
is $S$. The graph $[G:S]$ will play a role of a {\em state} in our
constructions. Let $b_0 ([G:S])$ denotes the number of connected
components of $[G:S]$ (that is the zeroth Betti number of the
graph), and $b_1([G:S])$ denotes the first Betti number of $[G:S]$.
The following state sum formula for the chromatic polynomial:
$$
P_G (\lambda) = \sum_{S \subseteq E(G)} \, (-1)^{|S|} \,
\lambda^{b_0([G:S])} = \sum_{i \geq 0} (-1)^i \sum_{S \subseteq
E(G), |S|=i} \, \lambda^{b_0([G:S])}
$$
was used in \cite{HR05} for its categorification. We denote the
chain complex constructed in \cite{HR05} by $\{C^i_P\}$.

The following well-known state sum formula for the Tutte polynomial
(see, e.g. \cite{We}):
$$
T(G; x,y) = \sum_{S \subseteq E(G)} \, (x-1)^{-b_0 (G) + b_0([G:S])}
\, (y-1)^{b_1 ([G:S])}
$$
was used in \cite{JR05} for a categorification of a version of the
Tutte polynomial. We denote by $\{C^i_T\}$ the chain complex
constructed in
\cite{JR05}.

In the present paper we
categorify a multiple of the Yamada
polynomial $h(G; x,y)$ by constructing of a corresponding chain
complex $\{C^i_Y\}$.

\section{The Yamada polynomial}

The {\em Yamada polynomial} of $G$, denoted by $h(G; x, y)$, is
defined by the following formula \cite{Ya}:
\begin{equation} \label{hxy}
h(G; x,y) \, =\, \sum_{F\subseteq E} {(-x)^{-|F|} \, x^{b_0(G-F)} \,
y^{b_1(G-F)}},
\end{equation}
where $F$ ranges over the family of all subsets of $E=E(G)$, and
$|F|$ is the number of elements in $F$; $b_0$ and $b_1$ are the
Betti numbers in dimensions 0 and 1. In particular for the empty
graph $G = \emptyset$ we have $h(\emptyset)=1$.

Let $S$ be the complement of $F $ in $E$, i.~e. $S= E-F$. We denote
by $[G:S]$ the graph with vertex set $V(G)$ and edge set $S$. Then
the Yamada polynomial can be written as follows:
\begin{equation} \label{hxys}
h(G; x,y) \, =\, \sum_{S \subseteq E} {(-x)^{-|E|+|S|} \,
x^{b_0([G:S])} \, y^{b_1([G:S])}}.
\end{equation}
It is obvious that $h(G; x,y)$ is a 2-variable Laurent polynomial in
$x$ and $y$ with nonnegative degrees on $y$.

Let us define a polynomial $\widetilde{g}(G,x,y)$ by the formula
\begin{eqnarray} \label{gxys}
\widetilde{g}(G,x,y) & = & (-x)^{|E|} h(G; x,y) \, \cr & = & \,
\sum_{S \subseteq E} (-1)^{|S|}{x^{|S|+ {b_0([G:S])}} \,
y^{b_1([G:S])}}.
\end{eqnarray}
Clearly, each monomial of $\widetilde{g}(G,x,y)$ has nonnegative
degrees on $x$ and $y$. Let us make change of variables $x=1+t$,
$y=1+w$ and define
\[
g(t,w)= \widetilde{g} (1+t, 1+w)= \sum_{S \subseteq E} (-1)^{|S|}
{(1+t)^{|S|+ {b_0([G:S])}} \, (1+w)^{b_1([G:S])}}.
\]

We intend to construct the chain complex and homology (in the sense
of Khovanov) corresponding to this polynomial.

The following evident statement was pointed out in \cite{JR05}.

\begin{lemma} \label{adding edge}
Let $G=(V,E)$ and $S$ be a subset of $E$. Suppose that $e \in E -
S$, and denote $S' = S \cup \{e\}$. Then
one of the following two cases occurs. \\
(i) If endpoints of $e$ belong to one component of $[G:S]$ then
$$
b_0([G : S'])  =  b_0([G:S]) \quad \text{and} \quad b_1([G : S']) =
b_1([G:S]) + 1.
$$
(ii) If endpoints of $e$ belong to different components of $[G:S]$
then
$$
b_0([G : S'])  =  b_0([G:S]) - 1 \quad \text{and} \quad  b_1([G :
S']) = b_1([G:S]).
$$
\end{lemma}
%\hfill $\square$

\section{The chain complex}

\subsection{Algebraic prerequisite}
 Let $R$ be a commutative ring with unit.
Recall (see \cite{CE} or \cite {ML} for example) that a $\Z$-{\em
graded} $R$-{\em module} or simply {\em graded} $R$-{\em module} $M$
is an $R$-module with a family of submodules $M_n$ such that $M$  is
a direct sum  $M=\oplus_{n\in \Z}M_n$. Elements of $M_n$ are called
{\em homogeneous elements of degree $n$}.

If $R=\Z$ and  $M=\oplus_{n\in \Z}M_n$ is a graded
$\mathbb{Z}$-module (abelian group) then the {\it graded dimension} of $M$ is the
power series
$$
q\dim M=\sum_n q^n \cdot \dim_{\mathbb{Q}} (M_n
\otimes \mathbb{Q}).
$$

In the same way a $\Z\oplus \Z$-{\it graded} or  {\it bigraded}
$R$-module is a $R$-module $M$ with a family of submodules
$M_{n,k}$, $n,k \in \Z$ such that
$$M=\oplus_{(n,k)\in \Z \oplus \Z }M_{n,k}.$$
Elements of $M_{n,k}$ are called {\em homogeneous elements of
bidegree $(n,k)$}. The {\em graded dimension} of $M$ over $\Z$ is
the 2-variable power series
$$
q\dim M=\sum_{n,k} x^n y^k \cdot \dim_{\mathbb{Q}} (M_{n,k} \otimes
\mathbb{Q}).
$$

\subsection{The general construction}

Let $M $ be a bigraded module over the ring $R$ equipped with an
associative multiplication
$$
m: M\otimes M \to M
$$
and a map
$$
u: R\to M,
$$
which even not necessary to be a unit for the multiplication $m$.
Let $N$ be any bigraded module over $R$. For each integer $\nu \geq
0$, let $$f_{\nu}: N^{\otimes \nu} \rightarrow N^{\otimes (\nu+1)}$$
be a degree preserving module homomorphism. Given such $M, N$ and
$f_{\nu}$, we can construct cohomology groups in the following
manner which is standard for Khovanov's approach.

As in Section~\ref{sec:pol_gr} we consider a graph $G = (V(G),
E(G))$ and $|E(G)|=n$. In Khovanov construction for links an
ordering of all crossings was done. Such an ordering is usual in
homological constructions. Here for graphs an ordering of edges of
$G$ is fixed:  $e_1, \cdots, e_n$. To visualise Khovanov
construction Bar-Natan \cite{BN05} suggests to consider the
$n$-dimensional cube with vertices $\{0,1\}^n$. For each vertex
$\alpha = (\alpha_1, \dots, \alpha_n)$ of the cube there corresponds
a subset $S=S_{\alpha}$ of $E(G)$, where $e_i \in S_{\alpha}$ if
$\alpha_i=1$. Bar-Natan  defines a {\it height} of the vertex
$\alpha = (\alpha_1, \dots, \alpha_n)$ as $|\alpha|=\sum \alpha_i$,
which is equal to the number of edges in $S_{\alpha}$.
 Each edge $\xi$ of the cube
$\{ 0, 1\}^E$ Bar-Natan  labels by a sequence $(\xi_1, \dots,
\xi_n)$ in $\{ 0, 1, *\}^E$ with exactly one ``$*$''. The tail
$\alpha_{\xi}(0)$ of $\xi$ is obtained by setting $*=0$ and the head
$\alpha_{\xi}(1)$ is obtained by setting $*=1$. The height $|\xi|$
is defined to be  equal to the number of 1's in the sequence
presenting $\xi$. We consider a subgraph $[G:S]$ of $G$ (see
Section~\ref{sec:pol_gr}) and take a copy of the $R$-module $M$ for
each edge of $S$ and  each connected component of $[G:S]$ and then
take a tensor product of copies of $M$ over the edges and the
components. Let $M^{\alpha}(G)$ be the resulting bigraded
$R$-module, with the bigrading induced from $M$. Thus,
$$
M^{\alpha}(G) \cong M^{\otimes \lambda}\otimes M^{\otimes \mu} ,
$$
where $\lambda=|S|$ and $\mu=b_0([G:S])$. Suppose
$N^{\alpha}(G)=N^{\otimes \nu}$, where $\nu = b_1([G:S])$. We define
$$
C^{\alpha}(G) = M^{\alpha}(G)\otimes N^{\alpha}(G) = M^{\otimes
\lambda} \otimes M^{\otimes \mu} \otimes N^{\otimes \nu}.
$$
So for each vertex $\alpha$ of the cube, we associated the bigraded
$R$-module $C^{\alpha}(G)$ (also denoted by $C^S(G)$, where
$S=S_{\alpha}$).
The {\em $i^{\text{th}}$ chain module} of the complex is
defined by
\begin{equation}
C^i(G):=\oplus_{|\alpha|=i} \, C^{\alpha}(G) .
\label{eq:Ci}
\end{equation}
The {\em differential maps}
$$
d^i: C^i(G)\rightarrow C^{i+1}(G)
$$
are defined using the multiplication $m$ on $M$, the map $u$, and
the homomorphisms $f_{\nu}$ as follows.

Consider the edge $\xi$ of the cube which joins two vertices
$\alpha_{\xi}(0)$ (the starting point) and $\alpha_{\xi}(1)$ (its
terminal). Denote the corresponding subsets of $E(G)$ by $S_0 =
S_{\alpha_{\xi}(0)}$ and $S_1 = S_{\alpha_{\xi}(1)}$. The edge of
the graph $e \in E(G)$ is such that $S_1 = S_0 \cup \{ e\}$. Let us
define now the {\em per-edge map}
$$
d_{\xi}: C^{\alpha_{\xi}(0)}(G) \rightarrow C^{\alpha_{\xi}(1)}(G).
$$
Denote $\lambda_i=|S_i|$, $\mu_i=b_0([G : S_i])$, and $\nu_i =
b_1([G : S_i])$ for $i=0,1$. Then we present
$$
d_{\xi} = d_{\xi}^M \otimes d_{\xi}^N : M^{\otimes \lambda_0}
\otimes M^{\otimes \mu_0} \otimes N^{\otimes \nu_0} \to M^{\otimes
\lambda_1} \otimes M^{\otimes \mu_1} \otimes N^{\otimes \nu_1} ,
$$
with $d_{\xi}^M : M^{\otimes \lambda_0} \otimes M^{\otimes \mu_0}
\to M^{\otimes \lambda_1} \otimes M^{\otimes \mu_1}$ and $d_{\xi}^N
: N^{\otimes \nu_0} \to N^{\otimes \nu_1}$.

Obviously, $\lambda_1 = \lambda_0 + 1$. Suppose that $d_{\xi}^M$
acts on the factor $M^{\otimes \lambda_0}$ of the tensor product
$M^{\otimes \lambda_0} \otimes M^{\otimes \mu_1}$ by the map $u$:
\begin{eqnarray*}
M^{\otimes \lambda_0} & = & M\otimes \dots \otimes M\otimes R
\otimes M \otimes \dots \otimes M \cr & \to & M\otimes \dots \otimes
M\otimes M \otimes M \otimes \dots \otimes M = M^{\otimes
(\lambda_0+1)} = M^{\otimes \lambda_1},
\end{eqnarray*}
where the position of $R$ is determined by the number of the
edge~$e$.

There are two possibilities which correspond to two cases in
Lemma~\ref{adding edge}.

If endpoints of $e$ belong to one component of $[G : S_0]$ then
$\mu_1 = \mu_0$ and $\nu_1 = \nu_0 + 1$. So, we put that $d_{\xi}^M$
acts on the factor $M^{\otimes \mu_0}$ of $M^{\otimes \lambda_0}
\otimes M^{\otimes \mu_0}$ by the identity map and $d_{\xi}^N:
N^{\alpha_{\xi}(0)}(G)  = N^{\otimes \nu_0} \to
N^{\alpha_{\xi}(1)}(G) = N^{\otimes \nu_1}$ acts by the homomorphism
$f_{\nu_0} : N^{\otimes \nu_0} \rightarrow N^{\otimes (\nu_0+1)}$.
Thus, the per-edge map $d_{\xi} = d_{\xi}^M \otimes d_{\xi}^N :
C^{\alpha_{\xi}(0)}(G) \to C^{\alpha_{\xi}(1)}(G)$ is defined.

If endpoints of $e$ belong to different components of $[G : S_0]$,
say $E_0$ and $E_1$, then $\mu_1 = \mu_0 -1$ and $\nu_1 = \nu_0$. In
this case we suppose that $d^M_{\xi}$ acts on the factor $M^{\otimes
\mu_0}$ of $M^{\otimes \lambda_0} \otimes M^{\otimes \mu_0}$ by the
multiplication map $m : M \otimes M \to M$ on tensor factors
corresponding to $E_0$ and $E_1$, and by the identity map on tensor
factors corresponding to remaining components. Put that $d_{\xi}^N:
N^{\alpha_{\xi}(0)}(G)  = N^{\otimes \nu_0} \to
N^{\alpha_{\xi}(1)}(G) = N^{\otimes \nu_1} = N^{\otimes \nu_0}$ acts
by the identity map. Thus, the per-edge map $d_{\xi} = d_{\xi}^M
\otimes d_{\xi}^N : C^{\alpha_{\xi}(0)}(G) \to
C^{\alpha_{\xi}(1)}(G)$ is defined.

Now we define the differential
$$
d^i:C^i(G) \rightarrow C^{i+1}(G)
$$
as usual by
$$
d^i = \sum_{|\xi|=i} \text{sign} (\xi) \,  d_{\xi} ,
$$
where $\text{sign} (\xi) =(-1)^{\sum_{i<j} \xi_i} $ and $j$ is the
position of ``$*$'' in the sequence \linebreak $(\xi_1, \dots ,
\xi_n)$ presenting $\xi$.

If $\Gamma$ is a subgraph of $G$ then there exists a chain
projection map
$$
p^i: C^i (G) \to C^i (\Gamma)
$$
defined by
$$ p^i(x) =
\begin{cases}
x ,
  &   \text{if $x\in C^\alpha$ such that  $S_\alpha \subset E(\Gamma)$ }, \\
0 ,
  & \text{otherwise. } \\
\end{cases}
$$

Denote the complex that we constructed by $\{ C^i_Y\}$.

The difference of our construction with that  of \cite{JR05} for $\{
C^i_T\}$ is the presence of the factor $M^{\otimes
\lambda}$ in each term $C^\alpha $ of $C^i$. We define a chain
map
$$
 \phi: C^i_T \to C^i_Y
$$
using the maps $u^{\otimes \lambda}$ on each term.

Suppose now that there exists a map
$$
\eta: M\to R
$$
such that its composition with $u$,
$$
\eta \circ u : R \to M\to R ,
$$
is identity. Then there exists a chain map
$$
\psi: C^i_Y \to C^i_T
$$
constructed using the maps $\eta^{\otimes \lambda}$ on each term.
The composition of $ \phi $ and $\psi$ is the identity map of $\{
C^i_T\}$ and so it becomes a direct summand of $\{ C^i_Y\}$. Denote
by $H^i_T(G)$ the cohomology theory constructed in \cite{HR05} and
by $H^i_Y(G)$ the cohomology theory defined by our
complex~$\{C^i_Y\}$.

\begin{theorem}\label{more general chain complex}
(a) The modules $C^i_Y(G)$ and the homomorphism $d^i$ form a chain
complex of bigraded modules whose differential preserves the
bidegree
$$0 \rightarrow C^0_Y(G) \stackrel{d^0}{\rightarrow}C^1_Y(G)
\stackrel{d^1}{\rightarrow} \cdots \stackrel{d^{n-1}}{\rightarrow}
C^n_Y(G) \rightarrow 0.$$
Denote it by
$C_Y(G)=C_{Y, M, N, f_{\nu}}(G)$. \\
(b) The cohomology groups $H^i_Y(G) (=H^i_{Y, M, N, f_{\nu}}(G))$
are invariants of the graph $G$, they are independent of
the ordering of the edges of $G$.  The isomorphism type of
the graded chain complex $C_Y(G)$ is an
invariant of $G$.\\
(c) If the graded dimensions of the modules $M$ and $N$ are well
defined, then the graded Euler characteristic is equal
\begin{eqnarray*}
\chi_{q}(C_Y(G)) & = & \mathrel{\mathop{\sum }\limits_{0\leq i\leq
n}}(-1)^{i} q\dim (H^{i}_Y) \\ & = & \mathrel{\mathop{\sum
}\limits_{0\leq i\leq n}} (-1)^{i} q\dim (C^{i}_Y) \\ & = &
g(G;q\dim M -1,q\dim N -1)
\end{eqnarray*}
(d) There is a morphism $\phi$ of chain complexes $C_T(G)\to C_Y(G)$
which generates a morphism of
graded modules $H^i_T(G) \to H^i_Y(G)$. \\
(e) If there exists a map $\eta: M\to R$ such that its composition
with $u$ is identity, then there exists a chain map
$$
 \psi: C^i_Y \to C^i_T
$$
such that its  composition with $\phi$ is the identity map of $\{
C^i_T\}$ and it becomes a direct summand of $\{ C^i_Y\}$. The same
is true for the cohomologies $H^i_Y(G)$ and $ H^i_T(G)$. \\
(f) The constructions above are functorial with respect to inclusions of
subgraphs $\Gamma \subset G$.
\end{theorem}

\begin{proof} We follow the proofs of analogous statements for
categorifications of the chromatic polynomial and the Tutte
polynomial from \cite{HR05} and \cite{JR05}.

(a) The map $d$ is degree preserving since it is built on the degree
preserving maps. It remains to show that $d \cdot d =0$. Let $S
\subseteq E(G)$. Consider the result of adding two edges $e_k$ and
$e_j$ to $S$ where $e_k$ is ordered before $e_j$. It is enough to
show that
\begin{equation}
d_{(...1...*...)} d_{(...*...0...)} = d_{(...*...1...)} d_{(...0...*...)}
\label{eq:com_dif}
\end{equation}

The proof of (\ref{eq:com_dif}) consists of checking various situations,
depending on how many components we have with or without $e_k$ and $e_j$.
Consider, for example the case when $e_k$ joins the edges of the same
component, and $e_j$ joins this component with the other one.
 Then we have
\begin{eqnarray*}
&  C^{S} (G)  =  M^{\otimes \lambda} \otimes M^{\otimes \mu} \otimes
N^{\otimes \nu} , & \\
& C^{S\cup  \{e_k\}} (G)  =  M^{\otimes (\lambda +1)} \otimes
M^{\otimes \mu} \otimes N^{\otimes (\nu+1)} , & \\
& C^{S \cup \{ e_j\}} (G)  =  M^{\otimes (\lambda + 1)} \otimes
M^{\otimes (\mu-1)} \otimes N^{\otimes \nu} , & \\
& C^{S \cup \{ e_k, e_j \}} (G)  =  M^{\otimes (\lambda +2)} \otimes
M^{\otimes (\mu - 1)} \otimes N^{\otimes (\nu+1)} , &
\end{eqnarray*}
and the per-edge maps act on factors of the tensor products as
follows:
\begin{eqnarray*}
d_{(...*...0...)} = (u, id, f_{\nu}),  \qquad
d_{(...1...*...)} = (u, m, id), \\
d_{(...0...*...)} = (u, m, id), \qquad
d_{(...*...1...)} = (u, id, f_{\nu}).
\end{eqnarray*}
This implies $d^{i} \cdot d^{i+1} =0$.

(b) The proof is the same as the proof of Theorem~2.12 in
\cite{HR05}. For any permutation $ \sigma $ of $\{1,..,n\},$
we define $G_{\sigma }$ to be the same graph but
with labels of edges permuted according to $ \sigma$. It is enough
to prove the result
when $\sigma =(k,k+1).$ Define an isomorphism $f$ of complexes
$$f: C^*(G)\to C^*(G_\sigma)$$
as follows.
 For any subset $S$ of $E$ with $i$ edges,
there is a summand in $C^i(G)$ and one in $C^i(G_{\sigma})$ that
defined by $S$.  Let $\alpha = (\alpha_1, \dots,
\alpha_n)$ be the vertex of the cube that corresponds $S$ in $G$ and
$f_S$ be the map between these two summands that is equal to $- id$ if
$ \alpha_{k} = \alpha_{k+1} = 1$ and equal to $id$ otherwise.
We define $f : C^i (G) \to C^i (G_{\sigma})$ d by $f = \oplus_{|S|=i} \,
f_{S}$. Obviously, $f$ is an isomorphism.

This shows that the isomorphism class of the chain complex
is an invariant of the graph.

(c) It follows from homological algebra that
$$
\sum_{0\leq i\leq n}(-1)^{i} q \dim (H^{i}(G))=\sum_{0\leq i\leq n}
(-1)^{i} q \dim (C^{i}(G)).
$$
We use (\ref{eq:Ci}) and the equality
$$
q \dim C^{\alpha}(G)
=(q\dim M)^{|S|+{b_0([G:S])}}(q\dim N)^{b_1([G:S])}
$$
which is exactly the contribution of the state $[G:S]$ in
$g(G;t,w)$.

Proofs of statements (d), (e), and (f) follow obviously from the
above considerations.
\end{proof}

\subsection{The special case}
Let $R=\Z$, and the role of $M$ and $N$ play the algebras
$A=\mathbb{Z}[t]/(t^2)$ and
$B=\Z[w]/(w^2)$, where $\deg
t=(1,0)$ and $\deg w=(0,1)$. Algebras $A$ and $B$ are bigraded
algebras with $ q\dim A=1+t$ and $q\dim B=1+w$. The map
$u_A: \Z\to A$ is given by
\begin{equation} \label{uA}
u_A(1) = 1,
\end{equation}
and the map $\eta_A: A\to \Z$ is given by
\begin{equation} \label{eA}
\eta_A(1) = 1, \quad \eta_A(t)= 0.
\end{equation}

Note that $A^{\otimes m}\otimes B^{\otimes n}$ is a bigraded
$\mathbb{Z}$-module whose graded dimension is $q\dim A^{\otimes
m}\otimes B^{\otimes n}=(1+t)^m (1+w)^n$. The algebra structure on
$B$ is not used, and the map
$$
B^{\otimes k}\to B^{\otimes (k+1)}
$$
is constructed by the map $ u_B : \Z\to  B$ is analogous to $u_A$:
$$
u_B(1) = 1.
$$

Applying Theorem~\ref{more general chain complex} to this case we
get

\begin{theorem}\label{chain complex}  The analogues of items (a) and (b) of
Theorem~\ref{more general chain complex} hold. The item
(c) is precised in the following form:\\
($c^{\prime}$) The graded Euler characteristic is equal
\begin{equation*}
\chi_{q}(C_Y(G))  =  \mathrel{\mathop{\sum }\limits_{0\leq i\leq
n}}(-1)^{i} q\dim (H^{i}_Y)  =  \mathrel{\mathop{\sum
}\limits_{0\leq i\leq n}} (-1)^{i} q\dim (C^{i}_Y) = {g}(G;t,w) .
\end{equation*}
As for the items (d) and (e) we have the following \\
($d^{\prime} \cup e^{\prime}$) There are morphisms  of chain
complexes $\phi: C_T(G)\to C_Y(G)$ and $\psi: C_Y (G) \to C_T (G)$
with the composition equals to the identity of $\{ C^i_T\}$, so it
becomes a direct summand of $\{ C^i_Y\}$. These morphisms generate
morphisms of graded modules $H^i_Y(G) \to H^i_T(G)$ and $H^i_Y(G)
\to H^i_T(G)$, with the composition equal to the identity of $\{
H^i_T\}$, so it becomes a direct summand of $\{ H^i_Y\}$.
\end{theorem}

\section{The Example}

Let us illustrate the above constructions for the graph $P_2$
consisting of two vertices connected by two edges, that is, the
bigon: \includegraphics[width=0.3in]{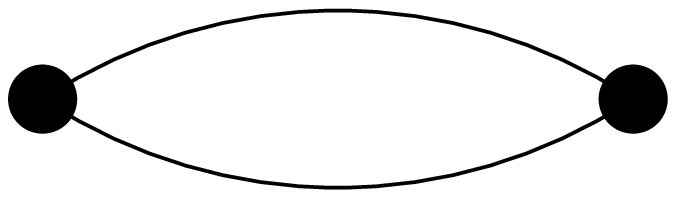}. Thus, $n=2$ and for vertices $\alpha$ of $\{0,1\}^2$
we get the following table:
$$
\begin{array}{|c|c|c|c|c|c|}
\hline \alpha & |\alpha| & \lambda & \mu & \nu & C^{\alpha} \\
\hline (0,0) & 0 & 0 & 2 & 0  & A \otimes A \\ \hline (1,0) & 1 & 1
& 1 & 0  & A \otimes A \\ \hline (0,1) & 1 & 1 & 1 & 0  & A \otimes
A \\ \hline (1,1) & 2 & 2 & 1 & 1  & A \otimes A \otimes A \otimes B \\
\hline
\end{array}
$$
Therefore, $C^0 = A \otimes A$, $C^1 = A \otimes A \oplus A \otimes
A$, $C^2 = A \otimes A \otimes A \otimes B$. The corresponding chain
complex is:
\begin{equation}\label{com_P2}
0 \rightarrow A \otimes A \overset{d^0}{\rightarrow} A \otimes A
\oplus A\otimes A \overset{d^1}{\rightarrow} A \otimes A \otimes
A\otimes B {\rightarrow} 0
\end{equation}
Where the differential map $d^0 = d_{(0,*)} + d_{(*,0)}$ acts as
follows:
$$
\begin{array}{rcl}
t\otimes t & \mapsto &  (0, 0) \\
1_A \otimes 1_A & \mapsto & (1_A\otimes 1_A, 1_A\otimes 1_A) \\
t \otimes 1_A & \mapsto & (t \otimes 1_A, t \otimes 1_A) \\
1_A \otimes t & \mapsto & (1_A \otimes t,1_A\otimes t)
\end{array}
$$
The kernel of $d^0$ is generated by the elements $t \otimes t$ and
$t \otimes 1_A- 1_A \otimes t$. Thus $H^0_Y (P_2) \cong A\{ (1, 0)\}
\cong \Z(1,0)\oplus \Z(2,0)$. Here $A\{ (1, 0)\}$ denotes the module
$A$ with the bidegrees shifted by $ (1, 0)$.

The differential map $d^1 = d_{(*,1)} - d_{(1,*)}$ acts as
following:
$$
\begin{array}{rcl}
(1_A\otimes 1_A,0)   & \mapsto &   -1_A\otimes 1_A\otimes 1_A
\otimes 1_B, \\
(1_A\otimes t,0)  &  \mapsto  & -1_A\otimes 1_A\otimes t \otimes 1_B, \\
(t\otimes 1_A,0)  & \mapsto  & -t\otimes 1_A\otimes 1_A \otimes 1_B,
\\ (t\otimes t,0)  & \mapsto &  -t\otimes 1_A\otimes t \otimes 1_B, \\
(0,1_A\otimes 1_A) & \mapsto & 1_A\otimes 1_A\otimes 1_A \otimes
1_B, \\ (0,1_A\otimes t)  & \mapsto &  1_A\otimes 1_A\otimes t \otimes 1_B, \\
(0, t\otimes 1_A) & \mapsto & 1_A\otimes t\otimes 1_A \otimes 1_B,
\\ (0,t\otimes t) & \mapsto & 1_A\otimes t\otimes t \otimes 1_B .
\end{array}
$$
The kernel of $d^1$ is generated by the elements $(1_A \otimes t,
1_A \otimes t)$ and $(1_A\otimes 1_A, 1_A\otimes 1_A)$. Two them lie
in the image of $d^0$, thus $H^1_Y (P_2) \cong 0$.  We have
$$H^2_Y(P_2) \cong \Z(2,0) \oplus \Z(3,0) \oplus \Z(0,1) \oplus
3\Z(1,1)\oplus 3\Z(2,1)\oplus \Z(3,1).$$ Clearly $H^i_Y(P_2)=0$, for
$i\geq 3$. Hence
\begin{eqnarray*}
\chi(H^*(P_2)) & = & t+ 2t^2+t^3+w+3tw +3t^2w+3t^3w \cr & = &
-(1+t)^2+(1+t)^3(1+w) \cr & = &  g(P_2;t,w) .
\end{eqnarray*}

Let us compare this with the cohomology of E.~F.~Jasso-Hernandez and
Y.~Rong \cite{JR05}. They have the complex
\begin{equation}
0 \rightarrow A \otimes A \overset{d^0}{\rightarrow} A
\oplus A \overset{d^1}{\rightarrow} A \otimes B \overset{d^2}{\rightarrow} 0
\end{equation}
which is evidently the direct composant of (\ref{com_P2}) as well as
cohomology groups Thus $H^0_T (P_2)
 \cong \Z(1,0)\oplus \Z(2,0)$, $H^1_T(P_2) \cong 0$,
$H^2_T(P_2) \cong  \Z(0,1) \oplus \Z(1,1)$.

\end{document}